\documentclass[10pt] {bcp}

\usepackage{amsfonts}

\theoremstyle{definition}

\begin{document}

\mathclass{Primary 41A36.}

\abbrevauthors{A.S. Serdyuk, T.A. Stepanyuk} \abbrevtitle{Lebesque-type inequalities }

\title{
 Lebesque-type inequalities  \\
for the Fourier sums on  classes of \\  generalized Poisson integrals 
 $ \ $}


\author{Anatoly S. Serdyuk}
\address{Institute of Mathematics NAS of
Ukraine \\
 Tereschenkivska st. 3,  01601   Kyiv, Ukraine\\
E-mail: serdyuk@imath.kiev.ua}

\author{Tetiana A. Stepanyuk}
\address{Graz University of Technology \\ Kopernikusgasse 24/II 8010, Graz, Austria\\
E-mail: tania$_{-}$stepaniuk@ukr.net}

\maketitlebcp

{\bf Summary}

\noindent For functions from the set of generalized Poisson integrals  $C^{\alpha,r}_{\beta}L_{p}$, $1\leq p <\infty$, we obtain upper estimates for the deviations of Fourier sums in the uniform metric in terms of the best approximations of the  generalized derivatives $f^{\alpha,r}_{\beta}$ of functions of this kind by trigonometric polynomials in the metric  of the spaces $L_{p}$. Obtained estimates are asymptotically best possible.

\vskip 5mm

Let $L_{p}$,
$1\leq p<\infty$, be the space of $2\pi$--periodic functions $f$ summable to the power $p$
on  $[0,2\pi)$, in which
the norm is given by the formula
$\|f\|_{p}=\Big(\int\limits_{0}^{2\pi}|f(t)|^{p}dt\Big)^{\frac{1}{p}}$; $L_{\infty}$ be the space of measurable and essentially bounded   $2\pi$--periodic functions  $f$ with the norm
$\|f\|_{\infty}=\mathop{\rm{ess}\sup}\limits_{t}|f(t)|$; $C$ be the space of continuous $2\pi$--periodic functions  $f$, in which the norm is specified by the equality
 ${\|f\|_{C}=\max\limits_{t}|f(t)|}$.

Denote by $C^{\alpha,r}_{\beta}L_{p}, \ \alpha>0, \ r>0, \ \beta\in\mathbb{R}, \ 1\leq p\leq\infty,$ the set of all  $2\pi$--periodic functions, representable for  all
$x\in\mathbb{R}$ as convolutions of the form (see, e.g., \cite[p.~133]{Stepanets1})
\begin{equation}\label{conv}
f(x)=\frac{a_{0}}{2}+\frac{1}{\pi}\int\limits_{-\pi}^{\pi}P_{\alpha,r,\beta}(x-t)\varphi(t)dt,
\ a_{0}\in\mathbb{R}, \  \ \varphi\perp1, 
\end{equation}
where $\varphi\in L_{p}$ and $P_{\alpha,r,\beta}(t)$ are fixed generated kernels

\begin{equation}\label{kernel}
P_{\alpha,r,\beta}(t)=\sum\limits_{k=1}^{\infty}e^{-\alpha k^{r}}\cos
\big(kt-\frac{\beta\pi}{2}\big), \ \  \ \alpha,r>0, \ \  \beta\in
    \mathbb{R}.
\end{equation}
The kernels  $P_{\alpha,r,\beta}$ of the form (\ref{kernel}) are called generalized Poisson kernels. For $r=1$ and $\beta=0$ the kernels $P_{\alpha,r,\beta}$ are usual Poisson kernels of harmonic functions.

If the functions $f$ and $\varphi$ are related by the equality (\ref{conv}), then function $f$ in this equality is called generalized Poisson integral of the function $\varphi$. The function $\varphi$ in equality (\ref{conv}) is called as generalixed derivative of the function $f$ and is denoted by $f^{\alpha,r}_{\beta}$.

The set of functions $f$ from $C^{\alpha,t}_{\beta}L_{p}$, $1\leq p\leq \infty$, such that 
 $ f^{\alpha,r}_{\beta}\in B_{p}^{0}$, where
$$
B_{p}^{0}=\left\{\varphi: \ ||\varphi||_{p}\leq 1, \  \varphi\perp1\right\},
 \ 1\leq p\leq \infty,
$$
we will denote by $C^{\alpha,r}_{\beta,p}$.

Let $E_{n}(f)_{L_{p}}$ be the best approximation of the function $f\in L_{p}$ in the metric of space $L_{p}$, $1\leq p\leq\infty$, by the trigonometric polynomials $t_{n-1}$ of degree $n-1$, i.e.,
$$
E_{n}(f)_{L_{p}}=\inf\limits_{t_{n-1}}\|f-t_{n-1}\|_{L_{p}}.
$$

Let  $\rho_{n}(f;x)$ be the following quantity
\begin{equation}\label{rhoF}
\rho_{n}(f;x):=f(x)-S_{n-1}(f;x),
\end{equation}
 where $S_{n-1}(f;\cdot)$ are the partial Fourier sums of order  $n-1$ for a function $f$.
 
 Least upper bounds of the quantity $\|\rho_{n}(f;\cdot)\|_{C}$  over the classes $C^{\alpha,r}_{\beta,p}$, we denote by $ {\cal E}_{n}(C^{\alpha,r}_{\beta,p})_{C}$, i.e.,
  \begin{equation}\label{sum}
 {\cal E}_{n}(C^{\alpha,r}_{\beta,p})_{C}=\sup\limits_{f\in
C^{\alpha,r}_{\beta,p}}\|f(\cdot)-S_{n-1}(f;\cdot)\|_{C},  \ r>0, \ \alpha>0, \ 1\leq p \leq \infty.
  \end{equation}
  
  Asymptotic behaviour of the quantities $ {\cal E}_{n}(C^{\alpha,r}_{\beta,p})_{C}$ of the form (\ref{sum}) was studied in \cite{Stepanets1984}--\cite{SerdyukStepanyuk2018}.
  
  In \cite{Stepanets1989N4}--\cite{SerdyukMusienko}  it was found the analogs of the Lebesque inequalities for functions $f\in C^{\alpha,r}_{\beta}L_{p}$ in the case $r\in(0,1)$ and $p=\infty$, and also in the case  $r\geq1$ and $1\leq p\leq\infty$, where the estimates for the deviations $\|f(\cdot)-S_{n-1}(f;\cdot)\|_{C}$ are expressed in terms of the best approximations $E_{n}(f^{\alpha,r}_{\beta})_{L_{p}}$. Namely, in \cite{Stepanets1989N4} it is proved that the following best possible inequalitiy holds
   \begin{equation}\label{Stepanets}
 \|f(\cdot)-S_{n-1}(f;\cdot)\|_{C}\leq \Big(\frac{4}{\pi^{2}}\ln n^{1-r}+O(1)\Big)e^{-\alpha n^{r}}E_{n}(f^{\alpha,r}_{\beta})_{L_{\infty}},
  \end{equation}
  where $O(1)$ is a quantity uniformly bounded with respect to $n$, $\beta$ and $f\in C^{\alpha,r}_{\beta}L_{\infty}$.
  
  The present paper is a continuation of \cite{Stepanets1989N4}, \cite{StepanetsSerdyuk}, and is devoted to getting asymptotically best possible analogs of Lebesque-type inequalities on the sets $C^{\alpha,r}_{\beta}L_{p}$, $r\in(0,1)$ and $p\in[1,\infty)$. This case was not considered yet. Let formulate the results of the paper.

By $F(a,b;c;d)$ we denote Gauss hypergeometric function
  $$
F(a,b;c;z)=1+\sum\limits_{k=1}^{\infty}\frac{(a)_{k}(b)_{k}}{(c)_{k}}\frac{z^{k}}{k!},
$$
$$
(x)_{k}:=\frac{x}{2}\Big(\frac{x}{2}+1\Big)\Big(\frac{x}{2}+2\Big)...\Big(\frac{x}{2}+k-1\Big).
$$

For arbitrary  $\alpha>0$, $r\in(0,1)$ and $1\leq p\leq\infty$ we denote by $n_0=n_0(\alpha,r,p)$ the smallest integer $n$ such that
\begin{equation}\label{n_p}
 \frac{1}{\alpha r}\frac{1}{n^{r}}+\frac{\alpha r \chi(p)}{n^{1-r}}\leq{\left\{\begin{array}{cc}
 \frac{1}{14},  & p=1, \\
\frac{1}{(3\pi)^3}\cdot\frac{p-1}{p}, & 1< p<\infty, \\
\frac{1}{(3\pi)^3}, & p=\infty, \
  \end{array} \right.}
\end{equation}
where $\chi(p)=p$ for $1 \leq p<\infty$ and $\chi(p)=1$ for $p=\infty$.

The following statement holds.

\textbf{Theorem 1.} \emph{Let $0<r<1$,  $\alpha>0$, $\beta\in\mathbb{R}$ and $n\in \mathbb{N}$. Then in the case $1< p<\infty$ for any function
 $f\in C^{\alpha,r}_{\beta}L_{p}$ and $n\geq n_0(\alpha,r,p)$, fthe following inequality is true:
$$
\|f(\cdot)-S_{n-1}(f;\cdot)\|_{C}\leq e^{-\alpha n^{r}}n^{\frac{1-r}{p}}\bigg(\frac{\|\cos t\|_{p'}}{\pi^{1+\frac{1}{p'}}(\alpha r)^{\frac{1}{p}}}F^{\frac{1}{p'}}\Big(\frac{1}{2}, \frac{3-p'}{2}; \frac{3}{2}; 1\Big)+
$$
\begin{equation}\label{Theorem1Ineq1}
+
\gamma_{n,p}\Big(\Big(1+\frac{(\alpha r)^{\frac{p'-1}{p}}}{p'-1}\Big) \frac{1}{n^{\frac{1-r}{p}}}+\frac{(p)^{\frac{1}{p'}}}{(\alpha r)^{1+\frac{1}{p}}}\frac{1}{n^{r}}\Big)\bigg)E_{n}(f^{\alpha,r}_{\beta})_{L_{p}}, \ \ \frac{1}{p}+\frac{1}{p'}=1,
\end{equation}
where $F(a,b;c;d)$ is Gauss hypergeometric function,
 and in the case $p=1$  for any function
 $f\in C^{\alpha,r}_{\beta}L_{1}$ and $n\geq n_0(\alpha,r,1)$, the following inequality is true:
 \begin{equation}\label{Theorem1Ineq2}
\|f(\cdot)-S_{n-1}(f;\cdot)\|_{C}\leq e^{-\alpha n^{r}}n^{1-r}\Big(
\frac{1}{\pi\alpha r}+\gamma_{n,1}\Big(\frac{1}{(\alpha r)^{2}}\frac{1}{n^{r}}+\frac{1}{n^{1-r}}\Big)\Big)E_{n}(f^{\alpha,r}_{\beta})_{L_{1}}.\ \
 \end{equation}
In (\ref{Theorem1Ineq1})  and (\ref{Theorem1Ineq2}), the quantity ${\gamma_{n,p}=\gamma_{n,p}(\alpha,r,\beta)}$ is such that ${|\gamma_{n,p}|\leq(14\pi)^{2}}$.}

\emph{Proof of Theorem 1.}
Let $f\in C^{\alpha,r}_{\beta}L_{p}$, $1\leq p\leq \infty$. Then, at every point $x\in \mathbb{R}$ the following integral representation is true:
\begin{equation}\label{repr}
\rho_{n}(f;x)=f(x)-S_{n-1}(f;x)=
\frac{1}{\pi}\int\limits_{-\pi}^{\pi}f^{\alpha,r}_{\beta}(t)P_{\alpha,r,\beta}^{(n)}(x-t))dt,
\end{equation}
  where
 \begin{equation}\label{kernelN}
P_{\alpha,r,\beta}^{(n)}(t):=
\sum\limits_{k=n}^{\infty}e^{-\alpha k^{r}}\cos\Big(kt-\frac{\beta\pi}{2}\Big),  \ 0<r<1, \ \alpha>0, \ \beta\in\mathbb{R}.
\end{equation}
 
The function $P_{\alpha,r,\beta}^{(n)}(t)$ is orthogonal to any trigonometric polynomial $t_{n-1}$ of degree not greater than $n-1$. Hence, for any polynomial $t_{n-1}$ from we obtain
\begin{equation}\label{for1}
\rho_{n}(f;x)=
\frac{1}{\pi}\int\limits_{-\pi}^{\pi}\delta_{n}(t)P_{\alpha,r,\beta}^{(n)}(x-t)dt,
\end{equation}
where
\begin{equation}\label{delta}
\delta_{n}(x)=\delta_{n}(\alpha,r,\beta,n;x):=f^{\alpha,r}_{\beta}(x)-t_{n-1}(x).
\end{equation}

Further we choose the polynomial $t_{n-1}^{*}$ of the best approximation of the function $f^{\alpha,r}_{\beta}$ in the space $L_{p}$, i.e., such that
$$
\| f^{\alpha,r}_{\beta}-t^{*}_{n-1}\|_{p}=E_{n}(f^{\alpha,r}_{\beta})_{L_{p}}, \ \ 1\leq p\leq\infty,
$$
to play role of $t_{n-1}$ in (\ref{for1}). Thus, by using the inequality
\begin{equation}\label{HolderIneq}
\bigg\|\int\limits_{-\pi}^{\pi}K(t-u)\varphi(u)du \bigg\|_{C}\leq \|K\|_{p'}\|\varphi\|_{p},
\end{equation}
$$
\varphi\in L_{p}, \ \ K\in L_{p'}, \ \ 1\leq p\leq\infty, \ \ \frac{1}{p}+\frac{1}{p'}=1
$$
(see, e.g., \cite[p. 43]{Korn}), we get
\begin{equation}\label{for2}
\|f(\cdot)-S_{n-1}(f;\cdot)\|_{C}\leq
\frac{1}{\pi}\|P_{\alpha,r,\beta}^{(n)}\|_{p'}E_{n}(f^{\alpha,r}_{\beta})_{L_{p}}.
\end{equation}

For arbitrary $\upsilon>0$ and ${1\leq s\leq \infty}$  assume
\begin{equation}\label{norm_j}
 \mathcal{I}_{s}(\upsilon):=\big\|\frac{1}{\sqrt{t^{2}+1}} \big\|_{L_{s}[0,\upsilon]},
\end{equation}
where
$$
\|f \|_{L_{s}[a,b]}=
   {\left\{\begin{array}{cc}
\bigg(\int\limits_{a}^{b}|f(t)|^{s}dt
\bigg)^{\frac{1}{s}}, & 1\leq s<\infty, \\
\mathop{\rm{ess}\sup}\limits_{t\in[a,b]}|f(t)|, \ & s=\infty. \
  \end{array} \right.}
$$

It follows from the paper  \cite{SerdyukStepanyuk2017}
 for arbitrary $r\in(0,1)$, $\alpha>0$, $\beta\in\mathbb{R}$, $1\leq s\leq\infty$, $\frac{1}{s}+\frac{1}{s'}=1$, $n\in\mathbb{N}$ and $n\geq n_0(\alpha,r,s')$
 the following  estimate holds
$$
\frac{1}{\pi}\|P_{\alpha,r,\beta}^{(n)} \|_{s}= \ \ \ \ \ \ \ \ \ \ \ \ \ \ \ \ \ \ \ \ \ \ \ \ \ \ \ \ \ \ \ \ \ \ \ \ \ \ \ \ \ \ \ \ \ \ \ \ \ \ \ \ \ \ \ \ \ \ \ \ \ \ \ \ \ \ \ \ \ \ \ \ \ \ \ \ \ \ \ \ \ \ \ \ \ \ \ \ \ \ \ \ \
$$
\begin{equation}\label{normKern}
=e^{-\alpha n^{r}}n^{\frac{1-r}{s'}}\bigg(\frac{\|\cos t\|_{s}}{\pi^{1+\frac{1}{s}}(\alpha r)^{\frac{1}{s'}}} \mathcal{I}_{s}\Big(\frac{ \pi n^{1-r}}{\alpha r}\Big)+
\delta_{n,s}^{(1)}\Big(\frac{1}{(\alpha r)^{1+\frac{1}{s'}}} \mathcal{I}_{s}\Big(\frac{ \pi n^{1-r}}{\alpha r}\Big)\frac{1}{n^{r}}+\frac{1}{n^{\frac{1-r}{s'}}}\Big)\bigg),
\end{equation}
where the quantity ${\delta_{n,s}^{(1)}=\delta_{n,s}^{(1)}(\alpha,r,\beta)}$,  satisfies the inequality  ${|\delta_{n,s}^{(1)}|\leq(14\pi)^{2}}$.

Substituting $s=p'=\infty$, from \ref{for2} and (\ref{normKern}) we get
(\ref{Theorem1Ineq2}).
 
 Further, according to \cite{SerdyukStepanyuk2017} for $n\geq n_{0}(\alpha,r,s')$,  $1<s<\infty$, $\frac{1}{s}+\frac{1}{s'}=1$, the following equality takes place
\begin{equation}\label{estimJ}
  \mathcal{I}_{s}\Big(\frac{ \pi n^{1-r}}{\alpha r}\Big)=  F^{\frac{1}{s}}\Big(\frac{1}{2}, \frac{3-s}{2}; \frac{3}{2}; 1\Big)+
  \frac{\Theta^{(1)}_{\alpha,r,s',n}}{s-1}\Big(\frac{\alpha r}{\pi n^{1-r}}\Big)^{s-1},
 \end{equation}
where $|\Theta^{(1)}_{\alpha,r,s',n}|<2$.

Let now consider the case $1<p<\infty$.

Formulas (\ref{normKern}) and (\ref{estimJ}) for $s=p'$ and $n\geq n_{0}(\alpha,r,p)$ imply
 $$
\frac{1}{\pi}\big\|\mathcal{P}_{\alpha,r,n} \big\|_{p'}= e^{-\alpha n^{r}}n^{\frac{1-r}{p}}\bigg(\frac{\|\cos t\|_{p'}}{\pi^{1+\frac{1}{p'}}(\alpha r)^{\frac{1}{p}}}F^{\frac{1}{p'}}\Big(\frac{1}{2}, \frac{3-p'}{2}; \frac{3}{2}; 1\Big)+
$$
$$
+
\gamma_{n,p}^{(1)}\Big(\frac{1}{p'-1}\frac{(\alpha r)^{\frac{p'-1}{p}}}{ n^{(1-r)(p'-1)}}+\frac{p^{\frac{1}{p'}}}{(\alpha r)^{1+\frac{1}{p}}}\frac{1}{n^{r}}+\frac{1}{n^{\frac{1-r}{p}}}\Big)\bigg)=
$$
$$
= e^{-\alpha n^{r}}n^{\frac{1-r}{p}}\bigg(\frac{\|\cos t\|_{p'}}{\pi^{1+\frac{1}{p'}}(\alpha r)^{\frac{1}{p}}}F^{\frac{1}{p'}}\Big(\frac{1}{2}, \frac{3-p'}{2}; \frac{3}{2}; 1\Big)+
$$
\begin{equation}\label{normKern1}
+
\gamma_{n,p}^{(2)}\Big(
\Big(1+\frac{(\alpha r)^{\frac{p'-1}{p}}}{p'-1}\Big)\frac{1}{ n^{(1-r)(p'-1)}}+\frac{p^{\frac{1}{p'}}}{(\alpha r)^{1+\frac{1}{p}}}\frac{1}{n^{r}}\Big)\bigg), 
\end{equation}
 where the quantities ${\delta_{n,p}^{(i)}=\delta_{n,p}^{(i)}(\alpha,r,\beta)}$,  satisfiy the inequality  ${|\delta_{n,p}^{(i)}|\leq(14\pi)^{2}}$, $i=1,2$.
 Estimate (\ref{Theorem1Ineq1}) follows from  (\ref{for2}) and (\ref{normKern1}).
 
 Theorem 1 is proved.
\vspace{-8mm}\begin{flushright} $\square$ \end{flushright}

It should be noticed, that estimates (\ref{Theorem1Ineq1}) and (\ref{Theorem1Ineq2}) are asymptotically best possible on the classes $C^{\alpha,r}_{\beta,p}$, $1\leq p<\infty$.

If $f\in C^{\alpha,r}_{\beta,p}$, then $\|f^{\alpha,r}_{\beta}\|_{p}\leq 1$, and 
$E_{n}(f^{\alpha,r}_{\beta})_{L_{p}}\leq 1$, $1\leq p<\infty$.
 Considering the least upper bounds of both sides of inequality (\ref{Theorem1Ineq1})  over the classes $C^{\alpha,r}_{\beta,p}$, $1<p<\infty$, we  arrive at the inequality
$$
{\cal E}_{n}(C^{\alpha,r}_{\beta,p})_{C}\leq e^{-\alpha n^{r}}n^{\frac{1-r}{p}}\bigg(\frac{\|\cos t\|_{p'}}{\pi^{1+\frac{1}{p'}}(\alpha r)^{\frac{1}{p}}}F^{\frac{1}{p'}}\Big(\frac{1}{2}, \frac{3-p'}{2}; \frac{3}{2}; 1\Big)+
$$
\begin{equation}\label{Theorem1Ineq1Remark}
+
\gamma_{n,p}\Big(\Big(1+\frac{(\alpha r)^{\frac{p'-1}{p}}}{p'-1}\Big) \frac{1}{n^{\frac{1-r}{p}}}+\frac{(p)^{\frac{1}{p'}}}{(\alpha r)^{1+\frac{1}{p}}}\frac{1}{n^{r}}\Big)\bigg)E_{n}(f^{\alpha,r}_{\beta})_{L_{p}}, \ \ \frac{1}{p}+\frac{1}{p'}=1.
\end{equation}

 Comparing this relation with the estimate of Theorem 4  from  \cite{SerdyukStepanyuk2017} (see also \cite{SerdyukStepanyuk2018}),
 we conclude that   inequality (\ref{Theorem1Ineq1}) on the classes $C^{\alpha,r}_{\beta,p}$, $1<p<\infty$, is asymptotically best possible.

  In the same way, the asymptotical sharpness of the estimate (\ref{Theorem1Ineq2}) on  the class $C^{\alpha,r}_{\beta,1}$ follows from comparing 
 inequality
\begin{equation}\label{Theorem1Ineq2Remark}
{\cal E}_{n}(C^{\alpha,r}_{\beta,p})_{C} \leq e^{-\alpha n^{r}} n^{1-r}\Big(
\frac{1}{\pi\alpha r}+\gamma_{n,1}\Big(\frac{1}{(\alpha r)^{2}}\frac{1}{n^{r}}+\frac{1}{n^{1-r}}\Big)\Big)E_{n}(f^{\alpha,r}_{\beta})_{L_{1}}
\end{equation}  
   and formula (9) from  \cite{SerdyukStepanyuk2016}.

\section*{Acknowledgements}

Second author is supported by the Austrian Science Fund FWF project F5503 (part of the Special Research Program (SFB)
"Quasi-Monte Carlo Methods: Theory and Applications")

\end{document}